\documentstyle[12pt]{article}
\title{\bf Boundedness
and Stability of Impulsively Perturbed
 Systems in a Banach Space}

\author
{L. Berezansky $^{1}$
 and  E. Braverman $^{2}$ }

\begin{document}
\maketitle

\footnotetext[1]{Ben-Gurion University of the Negev, Department of
Mathematics and Computer Science, Beer-Sheva 84105, Israel}
\footnotetext[2]{Technion-Israel Institute of Technology,
Department of Mathematics, Haifa 32000, Israel}

{\it Running head:}  Stability of impulsively perturbed systems

{\it Corresponding author:}  E. Braverman

{\it Correspondent's phone: } 972-429-4020

{\it Correspondent's e-mail:} maelena@tx.technion.ac.il

\begin{abstract}
We consider the influence of impulsive perturbations of a linear
impulsive equation in a Banach space on the existence of bounded
solutions and the exponential stability of the equation.
\end{abstract}

\section{Introduction}

{}~~~~~ The theory of impulsive differential equations
goes back to the work (Millman and Myshkis, 1960).
 These equations describe processes
changing their state abruptly at certain moments.
This means that the duration of the perturbation is negligible comparing with
the time of the process.
Perturbations of this kind occur in
control problems (Bressan and Rampazzo, 1991).
Recently the development of the theory of impulsive
differential equations in abstract spaces has begun
(Bainov {\it et al.}, 1988 both, 1989, 1993, Zabreiko {\it et al.},
1988).

We consider the problem
\begin{eqnarray}
 \dot{x} + A(t)x(t) = f(t), ~t \in [0, \infty),
\end{eqnarray}
\begin{eqnarray}
x(\tau_i +0) = B_i x(\tau_i - 0),
\end{eqnarray}
where $x(t)$ is an element of a Banach space $Y$,
$B_i: Y \rightarrow Y$ are linear bounded operators.
Our objective is to derive conditions such that the solution $x$ is
bounded on the half-line $[0, \infty)$ for any bounded
right-hand side $f$ :
\begin{eqnarray}
\sup_{t \geq 0} \parallel f(t) \parallel < \infty,
\end{eqnarray}
where $\parallel \cdot \parallel $ is the norm in $Y$.
In the control theory this problem is treated as follows.
Is the output $x$ bounded for any bounded input $f$ ?

Naturally the above problem is connected with
the exponential stability of a differential equation.
The work (Anokhin {\it et al.})
 deals with results of this kind for impulsive delay
differential equations.
The results for differential equations without impulses in
a Banach space can be found in (Dalecki\v{i} and Krein, 1974).
For impulsive differential equations in a Banach space a connection
of boundedness and stability is obtained in (Zabreiko
{\it et al.}, 1988).

Our approach is different.
We consider the homogeneous equation
\begin{eqnarray}
\dot{x}(t) + A(t)x(t) = 0, ~ t \in [0, \infty),
\end{eqnarray}
with non-homogeneous impulsive conditions
\begin{eqnarray}
x(\tau_i+0) = B_i x(\tau_i -0) + \alpha_i, ~i=1,2, \dots ,
\end{eqnarray}
where
\begin{eqnarray}
\sup_i \parallel \alpha_i \parallel < \infty.
\end{eqnarray}

Suppose any solution of (4),(5), with (6),
is bounded on the half-line.
We prove that under natural constraints each solution of (1),(2),
with the right-hand side satisfying (3), is also bounded on the half-line.

It is to be emphasized that unlike
(Zabreiko {\it et al.}, 1988)
we do not need to check whether the solution is bounded
for {\it any} bounded right-hand side.
Moreover, in the scalar case $Y= {\bf R}$ we can choose certain $\alpha_i$
(precisely, $\alpha_i~=~sign \{ \prod_{j=1}^i B_j \}$,
where $sign ~u = u/ \mid u \mid $ ).
Then by checking the boundedness of $x$ in the only case
we obtain that the solution of (1),(5) is bounded for any $f, \alpha$
satisfying (3),(6).
The result can be applied to differential equations without impulses
and it is new for them.

The paper is organized as follows.
First we obtain that if a solution of each problem (4),(5),
with (6), is bounded on the half-line,
then the evolution operator has an exponential estimate.
The section 3 deals with the main result described above.
In the section 4 the connection between impulsive equations and equations
without impulses is considered.

\section{Exponential estimates}

{}~~~~~Let $0 = \tau_0 < \tau_1 < \dots $ be fixed points,
$\lim _{i \rightarrow \infty} \tau _i = \infty , $
$Y$ be a Banach space with a norm $\parallel \cdot \parallel $,
by $\parallel \cdot \parallel $ we also denote the norm of a linear operator
acting in $Y$, $L(Y)$ is a space of linear bounded operators acting
in $Y$.

Let ${\bf l}_{\infty} (Y)$ be a space of sequences
$\alpha = \{ \alpha_i  \}_{i=1}^{\infty}, ~ \alpha_i \in Y, ~
i=1,2, \dots ,$ such that
$$ \sup_{i \geq 1} \parallel \alpha_i \parallel < \infty, $$
with the norm
$$ \parallel \alpha \parallel _{{\bf l}_{\infty} (Y)} =
\sup_{i \geq 1} \parallel \alpha_i \parallel . $$

$$ $$

\underline{\sl Definition.}
A function $x: {\bf R} \rightarrow Y$ absolutely continuous
in each $[\tau_i, \tau_{i+1}) $
is a solution of the impulsive equation (1),(5), if
for $t \neq \tau_i $ it satisfies (1) and for $t= \tau_i$
it satisfies (5).
\vspace{5 mm}

The solution is assumed to be right continuous.

Further we need the following hypotheses.

(H1) $ B_i: Y \rightarrow Y, ~ i=1,2, \dots $
are linear bounded operators ,
$A(t)$ is a continuous operator function with values in $L(Y)$,
$f(t)$ is a continuous function with values in $Y$.

(H2) $B_i : Y \rightarrow Y $ have bounded inverse operators.

(H3) There exists $\sigma > 0$ such that $\tau_{i+1} - \tau _i
< \sigma, ~ i=1,2, \dots ~$ .

(H4) There exists $\rho > 0$ such that $\tau_{i+1} - \tau _i
> \rho , ~ i=1,2, \dots ~$ .

Suppose (H1) and (H2) hold.
Then the impulsive equation (1),(5), $x(0)= \alpha$
 has one and only one solution
that can be presented as (Bainov, Kostadinov and Myshkis, 1988)
\begin{eqnarray}
x(t) = X(t) \alpha + \int_0^t C(t,s) f(s) ds +
\sum _{0< \tau_i \leq t} C(t, \tau_i) \alpha_i,
\end{eqnarray}
where the evolution operator
\begin{eqnarray}
C(t,s) = X(t) X^{-1} (s).
\end{eqnarray}

Here $X(t)$ is the solution of the operator equation
$$ \dot{X} (t) +A(t)X(t) = 0, $$
$$ X(\tau_i) = B_i X(\tau_i - 0), ~ i=1,2, \dots , $$
$$X(0)=I,$$
where $I : Y \rightarrow Y$ is the identity operator.

Obviously $C(t,s)$ satisfies the semi-group equality
\begin{eqnarray}
C(t,s) = C(t,\tau) C(\tau, s).
\end{eqnarray}

\newtheorem{guess}{Theorem}[section]
\begin{guess}
Suppose (H1),(H2) and (H3) hold and
 the solution of the problem (4),(5), $x(0) = 0$
is bounded for any $\alpha =
\{ \alpha_i \}_{i=1}^{\infty} \in {\bf l}_{\infty}(Y).$

Then there exist positive constants $N$ and $\nu$ such that the
inequality
\begin{eqnarray}
\parallel X(t) \parallel \leq N \exp (- \nu t)
\end{eqnarray}
 holds.
\end{guess}

{\sl Proof.}
By (7) the solution of the initial problem (4),(5), $x(0) = 0$
has the representation
\begin{eqnarray}
x(t) =
\sum _{0< \tau_i \leq t} C(t, \tau_i) \alpha_i .
\end{eqnarray}
Thus for any $t$ the right-hand side of (11) is a bounded
linear operator acting from ${\bf l}_{\infty} (Y) $
to $Y$ since
$$ \parallel x(t) \parallel \leq
\sum_{0 < \tau_i \leq t} \parallel C(t, \tau_i) \parallel
\cdot \parallel \alpha \parallel _{{\bf l}_{\infty} (Y) } , $$
where $ \alpha = (\alpha_1, \alpha_2, \dots , \alpha_i , \dots )
 \in {\bf l}_{\infty} (Y) $.
By the hypothesis of the theorem for any sequence
$\alpha \in {\bf l}_{\infty} (Y)$
the solution $x(t)$ is bounded.
Therefore the uniform boundedness principle
implies that there exists $k > 0$
such that
$$
\parallel x(t) \parallel
\leq k \parallel \alpha \parallel _{{\bf l}_{\infty}}
\mbox{~for any~} \alpha \in {\bf l}_{\infty} .
$$
Substituting of (8) in (11) gives
\begin{eqnarray}
 \parallel \sum _{0 < \tau_i \leq t}
X(t) X^{-1} (\tau_i) \alpha_i  \parallel
\leq k \parallel \alpha \parallel _{{\bf l}_{\infty}} .
\end{eqnarray}
Setting $\alpha_1 = X(\tau_1) y, \parallel y \parallel  =1 $
for $\tau_1 \leq t < \tau_2 $ we obtain
$$ \parallel  X(t) \parallel \leq k
\parallel  X(\tau_1) \parallel  $$
and $k \geq 1.$

By setting
$$ \alpha_1 = \frac{X(\tau_1) y}{\parallel X(\tau_1) \parallel },
\alpha_2 = \frac{X(\tau_2) y}{\parallel X(\tau_2) \parallel },
\parallel y \parallel = 1 , $$
we obtain from (12) for $t= \tau_2$
$$ \parallel X(\tau_2) y \parallel
\left[ (\parallel X(\tau_1) \parallel)^{-1} +
(\parallel X(\tau_2) \parallel)^{-1} \right] \leq k . $$
As $y$ is an arbitrary vector in ${\bf R}^n $ such that
$\parallel y \parallel = 1 $ then
$$ \parallel X(\tau_2) \parallel
\left[ (\parallel X(\tau_1) \parallel)^{-1} +
(\parallel X(\tau_2) \parallel)^{-1} \right] \leq k , $$
therefore
$$ 1 + \frac{\parallel X(\tau_2)\parallel}
{\parallel X(\tau_1)\parallel} \leq k. $$
Hence
\begin{eqnarray}
\parallel X(\tau_2)\parallel \leq  (k-1)
\parallel X(\tau_1)\parallel .
\end{eqnarray}
Without loss of generality we can assume $k > 1$.

By setting
$$ \alpha_1 = \frac{X(\tau_1) y}{\parallel X(\tau_1) \parallel },
\alpha_2 = \frac{X(\tau_2) y}{\parallel X(\tau_2) \parallel },
\alpha_3 = \frac{X(\tau_3) y}{\parallel X(\tau_3) \parallel },
\parallel y \parallel = 1 , $$
we obtain from (12) for $t = \tau_3$
$$ \parallel X(\tau_3) y \parallel
\left[ (\parallel X(\tau_1) \parallel)^{-1} +
(\parallel X(\tau_2) \parallel)^{-1} +
(\parallel X(\tau_3) \parallel)^{-1} \right] \leq k . $$
As $y \in {\bf R}^n$ is arbitrary such that
$\parallel y \parallel = 1$ then (13) gives
$$ 1 + \frac{\parallel X(\tau_3)\parallel}
{\parallel X(\tau_1)\parallel (k-1) }
 + \frac{\parallel X(\tau_3)\parallel}
{\parallel X(\tau_1)\parallel} \leq k. $$
Thus
$$ \parallel X(\tau_3)\parallel \leq
\parallel X(\tau_1)\parallel (k-1)^2 / k . $$

Now we prove by induction that
$$ \parallel X(\tau_{j+1})\parallel \leq
\parallel X(\tau_1)\parallel (k-1)^j / k^{j-1}  $$
for any positive integer $j$.

Suppose that
$ \parallel X(\tau_{i+1})\parallel \leq
\parallel X(\tau_1)\parallel (k-1)^i / k^{i-1}  , i \leq j .$

Then by setting in (12)
$$ \alpha_1 = \frac{X(\tau_1) y}{\parallel X(\tau_1) \parallel },
\dots  ,
\alpha_{j+1} = \frac{X(\tau_{j+1}) y}{\parallel X(\tau_{j+1}) \parallel },
\alpha_{j+2} = \frac{X(\tau_{j+2}) y}{\parallel X(\tau_{j+2}) \parallel },
$$ $\parallel y \parallel = 1 , $
we obtain for $t= \tau_{j+2}$
$$ 1 + \frac{\parallel X(\tau_{j+2}) \parallel}
{\parallel X(\tau_1) \parallel }
\left[ 1 + \frac{1}{k-1} + \frac{k}{(k-1)^2} +
\dots + \frac{k^{j-1}}{(k-1)^j} \right] \leq k. $$
As a sum of geometric progression
$$ \sum_{i=0}^{j-1} \frac{k^i}{(k-1)^{i+1} }
=\frac{\frac{k^j}{(k-1)^{j+1}} - \frac{1}{k-1}}
{\frac{k}{k-1} - 1} =
\frac{k^j}{(k-1)^j } - 1 , $$
then
$$ \frac{\parallel X(\tau_{j+2}) \parallel}
{\parallel X(\tau_1) \parallel }
\left[ 1 +  \frac{k^j}{(k-1)^j} - 1
\right] \leq k - 1 . $$
Hence $ \parallel X(\tau_{j+2}) \parallel
\leq  \parallel X(\tau_1) \parallel (k-1)^{j+1} / k^j . $

Let in (12) $$\alpha_1 = \alpha_2 = \dots = \alpha_j=0,
{}~\alpha_{j+1} = X(\tau_{j+1})y,
{}~\parallel y \parallel = 1. $$
Then
$$ \parallel X(t) \parallel \leq
k \parallel X(\tau_{j+1}) \parallel \leq
 \parallel X(\tau_1) \parallel (k-1)^j / k^{j-2} $$
for $\tau_{j+1} \leq t < \tau_{j+2}$.
Since by the hypothesis of the theorem $t<(j+2)\sigma,$
i.e. $j>t/{\sigma}-2$ ,  then
$$\ln ( \parallel X(t) \parallel ) \leq \ln (\parallel
X(\tau_1) \parallel ) - (t/{\sigma} - 2) \ln [k/(k-1)]
+ 2 \ln ~k.$$
By assuming
$$ \nu = \ln [k/(k-1)] /{\sigma}, $$
$$ N_1 = \parallel X(\tau_1) \parallel k^4 /(k-1)^2, $$
$$ N = \max \left\{ N_1,
\sup _{0 \leq t < \tau_1} [\exp (\nu t)
\parallel X(t) \parallel ] \right\} $$
we obtain the inequality (10) .
The proof of the theorem is complete.
$$ $$

\underline{\sl Definition.} The equation (1),(2)
is said to be exponentially stable
if there exist positive constants
$N$ and $\nu$ such that for any solution of the equation (4),(2)
the inequality
$$ \parallel x(t) \parallel \leq N \exp (- \nu t)
\parallel x(0) \parallel
 $$
holds.
$$ $$

\underline{\sl Corollary .}
Suppose the hypotheses of Theorem 2.1 hold.
Then the equation (1), (2) is exponentially stable.

$$ $$

Consider an additional hypothesis.

(H5) There exist positive constants $b$ and $M$ such that
\begin{eqnarray}
\sup_i \parallel B_i \parallel \leq b , ~~~~
\int_{\tau_i}^{\tau_{i+1}}
\parallel A(s) \parallel ds \leq M, ~ i=1,2, \dots ~.
\end{eqnarray}
$$ $$

\begin{guess}
Suppose the hypotheses (H1)-(H3) and (H5) hold and the solution
of the problem (4),(5), $x(0)=0$ is bounded for any
$\alpha = \{\alpha_i \}_{i=1}^{\infty} \in
{\bf l}_{\infty} (Y). $

Then there exist positive constants $N$ and $\nu$
such that
\begin{eqnarray}
\parallel C(t,s) \parallel \leq N \exp [- \nu (t-s)],
\end{eqnarray}
$0 \leq s < t < \infty $.
\end{guess}

{\sl Proof.}
We fix a positive integer $p$.
Then for any $u \in Y$ ~~$C(t, \tau_p ) u$ is the solution of the problem
(4),(5), $x(0)=0, ~ \alpha_i = 0, ~i \neq p, ~ \alpha _p = u $.

By repeating the proof of Theorem 2.1 we obtain
$$ \parallel C(t,\tau_p) \parallel \leq N_p \exp [- \nu (t-\tau_p)], $$
with
$$ N_p = \max \left\{
 \parallel C(\tau_{p+1} , \tau_p ) \parallel k^4 /(k-1)^2,
\sup _{\tau_p \leq t  < \tau_{p+1}  } [\exp (\nu (t- \tau_p ))
\parallel C(t, \tau_p ) \parallel ] \right\} . $$
Now we have to show that $N_p$ can be chosen
independently of $p$.

The operator $C(t, \tau_p)$ is the solution of the operator equation
$$ \dot{C} (t, \tau_p) + A(t) C(t, \tau_p) = 0 , ~
t \in [\tau_p, \tau_{p+1} ), $$
where $C(\tau_p , \tau_p) = I, ~ I:Y \rightarrow Y$ is the identity
operator,
$$ C(\tau_{p+1}, \tau_p) = B_{p+1} C(\tau_{p+1} - 0, \tau_p). $$

Thus for $t \in [\tau_p, \tau_{p+1})$ the operator $C(t, \tau_p)$
is the solution of the integral equation
$$ C(t, \tau_p) = I - \int_{\tau_p}^t A(\xi)C(\xi, \tau_p) d\xi .$$
By the Gronwall-Bellman inequality we obtain
$$ \parallel C(t, \tau_p) \parallel \leq
\exp \left\{ \int_{\tau_p}^t \parallel A(\xi) \parallel d \xi \right\}. $$
Therefore (14) gives
\begin{eqnarray}
\parallel C(t, \tau_p ) \parallel \leq e^M, ~~ t \in (\tau_p, \tau_{p+1}).
\end{eqnarray}
For $t= \tau_{p+1}$
\begin{eqnarray}
\parallel C(\tau_{p+1}, \tau_p) \parallel \leq
\parallel B_{p+1} \parallel \parallel C(\tau_{p+1} -0, \tau_p) \parallel
\leq b~e^M.
\end{eqnarray}
Hence
$$  \parallel C(t, \tau_p) \parallel \leq N_2 \exp [- \nu (t- \tau_p)], $$
where
$$ N_2 = \max \left\{ be^M k^4 /(k-1)^2 ,~ be^{\nu \sigma+M},~
e^{\nu \sigma +M} \right\}. $$

Now let $s \in [0, \infty)$ be arbitrary.
Suppose $\tau_p$ is the least of all $\tau_p \geq s$.
Then
$$ \parallel C(t,s) \parallel =
\parallel C(t,\tau_p) C(\tau_p, s) \parallel \leq
\parallel C(t,\tau_p) \parallel \parallel C(\tau_p, s) \parallel \leq $$
$$ \leq N_2 \exp [- \nu (t-\tau_p)] \cdot
N_2 \exp [- \nu (\tau_p - s)] =
N_2^2 \exp [- \nu (t-s)] . $$
The latter inequality may be proven similar to the proof of inequalities
(16),(17).

Therefore $\parallel C(t,s) \parallel \leq N \exp [- \nu (t-s)], $
with $N = N_2^2$.

The proof of the theorem is complete.

\section{Main results}

{}~~~~~Main results of this work are the following.

\begin{guess}
Suppose the hypotheses (H1)-(H5) hold and the solution of
(4),(5), $x(0)=0$ is bounded on the half-line for any
$\alpha = \{ \alpha_i \}_{i=1}^{\infty} \in {\bf l}_{\infty} (Y)$.

Then each solution $x(t)$ of the problem (1),(5) with $f, \alpha$
 satisfying (3),(6), is bounded for $t \geq 0$.
\end{guess}

{\sl Proof.}
By Theorems 2.1, 2.2 there exist positive constants $N$ and $\nu$
such that
$$ \parallel C(t,s) \parallel \leq N \exp [- \nu (t-s)], $$
$$ \parallel X(t) \parallel \leq N \exp (- \nu t). $$
The solution $x$ has the representation (7).
Thus
$$ \parallel x(t) \parallel \leq \parallel X(t) \parallel \parallel
x(0) \parallel + \sum_{0< \tau_i \leq t} \parallel C(t, \tau_i) \parallel
\parallel \alpha_i \parallel + $$
$$ + \int_0^t \parallel C(t,s) \parallel \parallel f(s) \parallel ds \leq
N \exp (- \nu t) \parallel x(0) \parallel + $$
$$ + \sum_{i=1}^{\infty} N \exp [- \nu (t- \tau_i)]
\parallel \alpha_i \parallel +
\int_0^t N \exp [- \nu (t-s)] \parallel f(s) \parallel ds. $$

Consider the second term.
Let $\tau_k$ be the greatest of all $\tau_i<t$.
Then
$$t- \tau_{k-1} \geq \tau_k-\tau_{k-1} > \rho, ~
t- \tau_{k-2} > 2 \rho, \dots, t- \tau_1 > (k-1) \rho. $$
Therefore
$$ \sum_{j=1}^{\infty} \exp [- \nu (t- \tau_j)] \leq
\sum_{j=0}^{\infty} \exp (- \nu j \rho) =
\frac{1}{1- \exp (- \nu \rho)} . $$
Hence
$$ \parallel x(t) \parallel \leq N \parallel x(0) \parallel +
\frac{N \parallel \alpha \parallel _{{\bf l}_{\infty} (X) }}
{1- \exp (- \nu \rho)} +
\frac{N}{\nu} \sup_{s \geq 0} \parallel f(s) \parallel .
$$
Thus $x(t)$ is bounded for $t \geq 0$, which completes the
proof of the theorem.
$$ $$

\begin{guess}
Suppose the hypotheses of Theorem 3.1 hold.

(a) Then
$$\lim_{n \rightarrow \infty} \parallel \alpha_n
\parallel = 0,
{}~\lim_{t \rightarrow \infty} \parallel f(t) \parallel = 0 $$
imply that for any solution $x$ of (1),(5)
$$ \lim_{t \rightarrow \infty} \parallel x(t) \parallel = 0 . $$

(b) Suppose that there exist positive constants $N_1$ and $\lambda$
such that
$$ \parallel \alpha_n \parallel \leq N_1 \exp (- \lambda n), ~
\parallel f(t) \parallel \leq N_1 \exp(- \lambda t). $$

Then for any solution $x$ of (1),(5) there exist positive
constants $N_0$ and $\nu_0$ such that
$$ \parallel x(t) \parallel \leq N_0 \exp(- \nu_0 t). $$
\end{guess}

{\sl Proof.}
(a) By Theorems 2.1, 2.2 as in the proof of Theorem 3.1
one obtains
$$ \parallel x(t) \parallel \leq
N \exp (- \nu t) \parallel x(0) \parallel +
\sum_{0< \tau_i \leq t} N \exp [- \nu (t-\tau_i)] \parallel \alpha_i
\parallel +
$$ $$
+ \int_0^t N \exp [- \nu (t-s)] \parallel f(s) \parallel ds . $$

We fix $ \varepsilon > 0$.
Let $t_1 > 0$ be such that
$$ t_1 > - \frac{1}{\nu} \ln \left\{
\frac{\varepsilon}{3N \parallel x(0) \parallel} \right\}. $$
Then the first term in the sum
$$ N \exp (- \nu t) \parallel x(0) \parallel < \varepsilon /3 .$$

Let $\tilde{t}$ be such that for $\tau_i \geq \tilde{t}$
$$ \parallel \alpha_i \parallel < \frac{\varepsilon}
{6N [1- \exp(- \nu \rho)]} $$
and $t_2 > \tilde{t}$ be such that for $t>t_2$
$$
\sum_{0< \tau_i \leq \tilde{t}} N \exp [- \nu (t
- \tau_i)] \parallel \alpha_i
\parallel  < \frac{\varepsilon}{6}. $$

Then similar to the proof of Theorem 3.1 we evaluate the second
term for $t>t_2$
$$\sum_{0< \tau_i \leq t} N \exp [- \nu (t-\tau_i)] \parallel \alpha_i
\parallel  =
\sum_{0< \tau_i \leq \tilde{t}} N \exp [- \nu (t-\tau_i)] \parallel \alpha_i
\parallel  + $$
$$+ \sum_{ \tilde{t} < \tau_i \leq t}
N \exp [- \nu (t-\tau_i)] \parallel \alpha_i
\parallel \leq
\frac{\varepsilon}{6}
+ \sum_{j=0}^{\infty} N \exp (- \nu j \rho)
\frac{\varepsilon}{6 N [1 - \exp( - \nu \rho)]}
< $$ $$ < \frac{\varepsilon}{6}
+ \frac{\varepsilon}{6}
= \frac{\varepsilon}{3} . $$

The corresponding estimation of the third term
$$\int_0^t N \exp[-\nu (t-s)]
\parallel f(s) \parallel ds
\mbox{ ~~~for~~~ }
\lim_{s \rightarrow \infty}
\parallel f(s) \parallel = 0 $$
is given in (Massera and Sch\"{a}ffer, 1966).
Thus for a certain $t_3$
$$\int_0^t N \exp[-\nu (t-s)]
\parallel f(s) \parallel ds < \frac{\varepsilon}{3},
$$
 for each $t>t_3$.

Therefore for $t> \max \{t_1,t_2,t_3\} $
$$ \parallel x(t) \parallel <
 \frac{\varepsilon}{3}
+ \frac{\varepsilon}{3}
+ \frac{\varepsilon}{3}
=  \varepsilon . $$
Since $\varepsilon > 0$ is arbitrary we obtain
 $ \lim_{t \rightarrow \infty} \parallel x(t) \parallel = 0$ .

(b) For obtaining an exponential estimate of $x(t)$
we have to evaluate each of three above terms.

The first term
$$ N \exp (- \nu t) \parallel x(0) \parallel
$$
 obviously has the exponential estimate,
the third one under the hypotheses of the corollary
also can be estimated as (Massera and Sch\"{a}ffer, 1966)
$$\int_0^t N \exp [- \nu (t-s)] \parallel f(s) \parallel
ds  \leq N_3 ~ \exp (- \nu_3 t), $$
with certain positive constants $N_3$ and $\nu_3$.

Now we prove that there exist positive constants $N_2$
and $\nu_2$ such that
$$ \sum_{0< \tau_i \leq t} N \exp [ - \nu (t- \tau_i)]
\parallel \alpha_i \parallel <
N_2 \exp (-\nu_2 t). $$
Let $\tau_k \leq t < \tau_{k+1}$.
Then
denoting
$$ \beta = \sum_{0< \tau_i \leq t}
N \exp [- \nu (t- \tau_i)] \parallel \alpha_i \parallel $$
we obtain
$$ \beta \leq \sum_{i=1}^k N \exp [- \nu (t- \tau_i)] N_1 \exp(-
\lambda i). $$

Since $t- \tau_i \geq (k-i) \rho $ ( see the proof of Theorem 3.1)
then
$$ \beta \leq N N_1 \sum_{i=1}^k \exp[- \nu (k-i) \rho ]
\exp (- \lambda i) = $$
$$ = N N_1 \exp (- \nu \rho k) \sum_{i=1}^k
\exp[( \nu \rho - \lambda) i] \leq $$
$$ \leq N N_1 \exp (- \nu \rho k) \frac{\exp [(\nu \rho - \lambda)
(k+1)]}{\exp ( \mid \nu \rho - \lambda \mid ) - 1} = $$
$$ = \frac{N N_1 \exp (\nu \rho)}
{\exp ( \mid \nu \rho - \lambda \mid ) - 1} \exp [- \lambda (k+1)]. $$
The hypothesis (H3) gives that $t< \tau_{k+1} < (k+1) \sigma$.
Therefore
$$ \beta \leq
 \frac{N N_1 \exp (\nu \rho)}
{\exp ( \mid \nu \rho - \lambda \mid ) - 1} \exp \left( -
\frac{\lambda}{\sigma} t  \right) = N_2 \exp(- \nu_2 t), $$
with $\nu_2 = - \lambda / \sigma ,$
$$
 N_2 = \frac{N N_1 \exp (\nu \rho)}
{\exp ( \mid \nu \rho - \lambda \mid ) - 1} . $$

Denoting
$$ N_0 = \max \{ N,N_2, N_3 \}, $$
$$ \nu_0 = \min \{ \nu ,\nu_2, \nu_3\} $$
we obtain
$$ \parallel x(t) \parallel \leq N_0 \exp (- \nu_0 t ), $$
which completes the proof.

$$ $$

\underline{\sl Corollary 3.3.}
Suppose $Y= {\bf R} $ (scalar case), (H1)-(H5) hold
 and the solution of (4),
\begin{eqnarray}
x(\tau_i) = B_i x(\tau_i - 0) + sign [B_i X(\tau_i - 0)], ~i=1,2, \dots,
\end{eqnarray}
$x(0)=0$
is bounded on the half-line.
Here $X$ is the function from the representation (7),
$~sign ~u = u/ \mid u \mid $.

Then each solution of (1),(5) is bounded
on $[0, \infty)$ for any $f, \alpha $ satisfying (3) and (6) correspondingly.

{\sl Proof.}
By the hypothesis of the corollary
there exists $Q > 0$ such that $\mid x(t) \mid \leq Q,
 ~t \geq 0,$
where $x$ is the solution of (4),(18), $x(0)=0$.

As $\alpha_i = X(\tau_i) / \mid X(\tau_i) \mid $,
then (7) implies
$$\mid x(t) \mid = \mid \sum_{0 < \tau_i \leq t}
C(t, \tau_i) X(\tau_i) / \mid X(\tau_i) \mid ~ \mid \leq Q. $$
However
$$ C(t, \tau_i) X(\tau_i) $$
have the same sign for any $i =1,2, \dots , \tau_i \leq t, $
for a fixed $t$ .
Therefore
$$ \mid \sum_{0 < \tau_i \leq t}
C(t, \tau_i) X(\tau_i) / \mid X(\tau_i) \mid ~ \mid =
\sum_{0 < \tau_i \leq t}
 \mid C(t, \tau_i) \mid  \leq Q. $$
Hence for any $\alpha = \{ \alpha_i\} _{i=1}^{\infty}
\in {\bf l}_{\infty} $
the solution $x$ of (4),(5) is bounded :
$$ \mid  x(t) \mid =
 \mid \sum_{0 < \tau_i \leq t}
 C(t, \tau_i)  \alpha _i \mid  \leq
\sum_{0 < \tau_i \leq t}
 \mid C(t, \tau_i) \mid
\cdot \parallel \alpha \parallel _{{\bf l}_{\infty}} \leq
Q \parallel \alpha \parallel _{{\bf l}_{\infty}} .
$$
Thus the hypotheses of Theorem 3.1 are satisfied.
Therefore each solution of (1),(5) is bounded on the half-line whenever
$f, \alpha$ satisfy (3),(6) correspondingly.

$$ $$

One can apply Theorems 2,1, 2.2 to ordinary differential
equations without impulses.
$$ $$

\underline{\sl Corollary 3.4.}
Suppose there exists a positive constant $\eta$
such that for any sequence
 $\alpha = \{ \alpha_i\} _{i=1}^{\infty}
\in {\bf l}_{\infty} (Y) $
the equation
$$\dot{x}(t) + A(t)x(t) = \sum_{i=1}^{\infty}
\alpha_i \delta(t- \eta i)
$$
has a bounded solution.

Then there exist positive constants $N$ and $\nu$
such that for the solution of (4) the estimates (10)
and (15) are valid.

Here $\delta (t -\eta i)$ is a delta function,
the derivative in the left-hand side and the equality
are understood in the distributional sense.
$$ $$

\underline{\sl Remark.}
Let $\eta (a,b)$ be the number of points $\tau_i$ lying
in the interval $[a,b)$.
It is to be noted that the usual condition (Bainov {\it et al.}, 1989)
$$ \lim _{\omega \rightarrow \infty}
\frac{\eta (t, t+ \omega)}{\omega} = q < \infty $$
is more restrictive than both hypotheses (H3) and (H4).

\section{Exponential estimates of impulsive and continuous
solutions}

{}~~~~~Since many results on exponential estimation of
the differential equation (1) are known the following
question is of special interest.
Let any solution of the problem (1), $x(0) = 0$
be bounded on $[0, \infty )$ for any bounded on $[0, \infty )$
right-hand side $f$.
Does this property preserve for the impulsive
problem (1),(2) ?
Boundedness of solutions is connected with exponential estimates
of an evolution operator.
The following examples illustrate that an exponential
estimate for the impulsive equation does not imply an exponential
estimate for the corresponding differential equation (1) and vice versa.

{\bf Example 1.}
The solution of the problem
$$ \dot{x} = 0, ~ x(i) = 0.5 x(i-0),~ i=1,2, \dots ,~ x(0)=1 $$
can be estimated
$$x(t) \leq 2 \exp (-t \ln 2), $$
while the solution of
$$ \dot{x} = 0, ~ x(0)=1 $$
is constant.

{\bf Example 2.}
The equation
$\dot{x} + x = f $
is exponentially stable.
However for the impulsive equation
$$\dot{x} + x = 0, ~x(i) = e~ x(i-0), ~i = 1,2, \dots , $$
we obtain $x(i) = x(0)$ .

It is well known that the operator $X(t)$ in the representation (7)
for the solution  of the impulsive equation
(1),(5) can be defined if we know the corresponding operator $U(t)$ of
the nonimpulsive equation (1)
(Bainov {\it et al.}, 1989)
\begin{eqnarray}
X(t) = U(t) U^{-1} (\tau_i)
\prod_{j=1}^i B_{i-j+1} X(\tau_{i-j+1}) X^{-1} (\tau_{i-j}),
\end{eqnarray}
for $t \in [\tau_i, \tau_{i+1} )$.
If we know the evolution operator $G(t,s) = U(t) U^{-1} (s) $
of the non-impulsive equation, then the evolution operator
of the impulsive equation
(Bainov {\it et al.}, 1989)
\begin{eqnarray}
C(t,s) = \left\{ \begin{array}{l}
G(t,s), ~ \tau_i \leq t,s < \tau_{i+1}, \\
G(t, \tau_i) \left[ \prod_{j=i}^{k+1}
B_j G(\tau_j, \tau_{j-1}) \right] B_k G(\tau_k, s) ,
\\ \tau_{k-1} \leq s < \tau_k < \tau_i \leq t < t_{i+1} , \\
G(t, \tau_i) \left[ \prod_{j=i}^{k-1}
B_j^{-1} G(\tau_j, \tau_{j+1}) \right] B_k^{-1} G(\tau_k, s) ,
\\ \tau_{i-1} \leq t < \tau_i < \tau_k \leq s < t_{k+1} .
\end{array}
\right.
\end{eqnarray}

$$ $$

\begin{guess}
 Let $X={\bf R}$ (scalar case).
Suppose there exists a positive constant
 $\varepsilon$ such that
$$ \mid \prod _{k=i}^j B_k \mid \geq \varepsilon $$
for each $i,j, ~i,j = 1,2, \dots $, (H1)-(H3) hold
and any solution of the initial problem
(4),(5), $x(0) = 0$ is bounded
on $[0, \infty )$
for any bounded sequence $\{ \alpha_i \} _{i=1}^{\infty}. $

Then there exist positive constants $N$ and $\nu$
such that the solution $U(t)$
of non-impulsive equation (4), $U(0)=1,$ satisfies (10).

If, in addition, (H5) holds, then for $G(t,s)$
the estimate (15) is valid for certain $N, \nu > 0$.
\end{guess}

{\sl Proof.}
By Theorem 2.1 the hypotheses of the theorem imply
that there exist $N_1 > 0, \nu > 0 $ such that
$$ \mid X(t) \mid \leq N_1 \exp (- \nu t), $$
where $X(t)$ is the solution of (4),(2), x(0)=1.
In the scalar case (19) can be rewritten as
\begin{eqnarray}
X(t) = U(t) \prod_{0 < \tau_i \leq t} B_i .
\end{eqnarray}
Therefore
$$ \mid U(t) \mid \leq \frac{\mid X(t) \mid }{\varepsilon}
\leq \frac{N_1}{\varepsilon} \exp (-\nu t) \leq N \exp (- \nu t) , $$
with $N = N_1 / \varepsilon $.

Similarly (20) implies the estimate (15) for $G(t,s)$.

$$ $$

\begin{guess}
Let the operators $B_i$ and $G(t,s)$ be commuting for each
$t,s~\geq~0, ~ i=1,2, \dots $.
Suppose any solution of the initial problem (1), $x(0) = 0$
is bounded on $[0, \infty )$ for any bounded on  $[0, \infty )$
right-hand side $f$ and there exists $Q>0$ such that
$$ \prod_{k=i}^j  \parallel B_k \parallel \leq Q $$
for each $i,j = 1,2, \dots .$

Then there exist positive constants $N$ and $\nu$ such that
the solution $X$ of the problem (4),(2), $X(0)=I,$  has an estimate
(10).

If, in addition, there exists $M > 0$ such that
$$ \int_i^{i+1} \parallel A(s) \parallel ds < M, ~ i = 1,2, \dots , $$
then for $C(t,s) = X(t) X^{-1} (s) $ the estimate (15)
is valid , with certain $N, \nu~>~0 $.
\end{guess}

{\sl Proof.}
The hypotheses of the theorem imply (Dalecki\v{i} and Krein,
1974, p.127)
that there exist $N_1, \nu > 0$ such that
$$ \parallel U(t) \parallel \leq N_1 \exp (- \nu t). $$
Since $B_i$ and $U(t)=G(t,0)$ are commuting,
we obtain (21).
Hence
$$ \parallel X(t) \parallel \leq \parallel U(t) \parallel
 \prod_{0 < \tau_i \leq t
} \parallel B_i \parallel \leq Q N_1 \exp (- \nu t) = N \exp (- \nu t), $$
with $N = N_1 Q .$

Similarly by applying (20) one obtains (15) ,
which completes the proof of the theorem.

\vspace{5 mm}

{\bf\large Acknowledgement}

\vspace{5 mm}

This work was supported by Israel Ministry of Science
and Technology and by the Centre for Absorption in Science,
Ministry of Immigrant Absorption State of Israel.
$$ $$

{\bf\large References}
\vspace{5 mm}

Anokhin, A., Berezansky, L., and Braverman, E.
Exponential stability of linear delay impulsive differential
equations . Preprint FUNCT-AN/9311004.

Bainov, D. D., Kostadinov, S. I., and Myshkis, A. D. (1988).
{\it Differential and Integral Equations, } {\bf 1} (2) ,
233.

Bainov, D.D., Kostadinov, S.I, Nguyen Van Minh, and Zabreiko, P.P. (1993).
{\it International Journal of Theoretical Physics}, {\bf 32}, 1275.

Bainov, D. D., Kostadinov, S. I., and Zabreiko, P. P. (1989).
{\it International Journal of Theoretical Physics,} {\bf 28}, 797.

Bainov, D. D., Zabreiko, P. P., and Kostadinov, S. I. (1988).
{\it International Journal of Theoretical Physics,} {\bf 27} , 373.

Bressan, A., and Rampazzo, F. (1991).
{\it Journal of Optimization Theory and Applications, } {\bf 71}
(1),
67.

Dalecki\v{i}, Ju. L., and Kre\v{i}n, M. G. (1974).
{\it Stability of Solutions of Differential Equations
in Banach Spaces,} AMS, Providence, Rhode Island.

Massera, J. L., and Sch\"{a}ffer, J. (1966).
{\it Linear Differential Equations and Function Spaces },
Academic Press, New York.

Zabreiko, P. P., Bainov, D. D., and Kostadinov, S. I. (1988).
{\it International Journal of Theoretical Physics,} {\bf 27}, 731.

\end{document}